\documentclass[12pt,a4paper]{article}

\usepackage{amssymb}

\usepackage{amsmath}

\newtheorem{theorem}{Theorem}[section]
\newtheorem{definition}[theorem]{Definition}
\newtheorem{lemma}[theorem]{Lemma}
\newtheorem{corollary}[theorem]{Corollary}

\newtheorem{example}[theorem]{Example}
\begin{document}
\title{Gr\"{o}bner-Shirshov Bases for Lie Algebras: after A. I. Shirshov }

\author{{\small L. A. Bokut\footnote {Supported by the RFBR and
the Integration Grant of the SB RAS (No. 1.9).}}\\
{\small School of Mathematical Sciences, South China Normal
University,}\\
{\small Guangzhou 510631, China}\\
{\small Sobolev Institute of Mathematics, Russian Academy of
Sciences,}\\
{\small  Siberian Branch, Novosibirsk 630090, Russia}\\
{\small bokut@math.nsc.ru }\\
\\
 {\small Yuqun Chen \footnote
{Corresponding author.} \footnote{Supported by the National Natural
Science Foundation of China (Grant No.10471045; 10771077) and the
Natural Science
Foundation of Guangdong Province (Grant No.021073; 06025062).}}\\
{\small  School of Mathematical Sciences, South China Normal
University,}\\
{\small Guangzhou 510631, China}\\
 {\small  yqchen@scnu.edu.cn}}

\date{}
 \maketitle

\maketitle \noindent\textbf{Abstract:} In this paper, we review
Shirshov's method for free Lie algebras invented by him in 1962
\cite{Sh} which is now called the Gr\"{o}bner-Shirshov bases theory.

\maketitle \noindent\textbf{Key words: } Lie algebra;
Lyndon-Shirshov word; Gr\"{o}bner-Shirshov basis.


\section{Introduction}
What is now called the Gr\"{o}bner-Shirshov method for Lie
 algebras invented by A. I. Shirshov in 1962 \cite{Sh}. Actually, that
 paper based on his paper \cite{S58} when Shirshov invented a new
 linear basis for a free Lie algebra which is now called Lyndon-Shirshov
 basis (it was defined independently in the paper \cite{C} in the same year).
We remark that Lyndon--Shirshov basis is a particular case of a
 series of bases of a free Lie algebra invented by A. I. Shirshov in his
 Candidate Science Thesis (Moscow State University, 1953, and his adviser was A. G. Kurosh)
 and published in 1962 \cite{S62-}
(cf. \cite{R} where these bases are called Hall Bases).
 We now cite the Zbl review by P. M. Cohn \cite{Cohn} of the paper \cite{S58}:
  ``The author varies the usual construction of basis commutators
 in Lie rings by ordering words lexicographically and not by length.
 This is used to give a very short proof of the theorem
 (Magnus \cite{M}, Witt \cite{W})
 that the Lie algebra obtained from a free associative algebra is
 free. Secondly he derives Friedrich's criterion (this Zbl 52,45)
 for Lie elements. As the third application he proves that every Lie
 algebra $L$ can be embedded in a Lie algebra $M$ such that in $M$
 any subalgebra of countable dimension is contained in a 2-generated
 subalgebra." We would like to add that it was also a beginning
 of Gr\"{o}bner-Shirshov bases theory for Lie and associative
 algebras. Lemma 4 of the paper, on special bracketing of a regular
 (Lyndon-Shirshov) associative word with a fix regular subword,
 leads to the algorithm of elimination of the leading
word of one Lie polynomial in other Lie polynomial, i.e., to the
reduction procedure, which is very familiar in the cases of
associative and associative-commutative polynomials. Also the above
Lemma 4 leads to the crucial notion of composition of two Lie
polynomials that will be defined lately in \cite{Sh}.

   \ \  As for paper \cite{Sh} itself, it is a fully pioneer paper in the
 subject. He defines a notion of the composition $(f,g)_w$ of two
 Lie (associative) polynomials relative to an associative word $w$
 (it was called lately by $S$-polynomial for commutative polynomials by
 B. Buchberger \cite{bu65} and \cite{bu70}).
 It leads to the algorithm for the construction of a Gr\"{o}bner-Shirshov
 basis ($GSB(S)$) of Lie (associative) ideal generated by some set
 $S$: to joint to S all nontrivial compositions and to eliminate the
leading monomials of one polynomial of S in  others. Shirshov proved
the lemma, now known as the Composition, or
 Composition-Diamond Lemma, that if $f \in Id_{Lie}(S)$, then
 $\overline{f}$, the leading associative word of $f$, has a form $\overline{f}=u
 \overline{s}v$, where $s \in GSB(S), \ u,v \in X^*$. Several years
 later, Bokut formulated this lemma in the modern form (see
\cite{b72}). Let $S$ be a set of Lie polynomials that is complete
under
 composition (i.e., any composition of polynomials of $S$ is trivial;
  on the other word, $S$ is a Gr\"{o}bner-Shirshov
 basis). Then if $f \in Id_{Lie}(S)$, then $\overline{f}=u
 \overline{s}v$, where $s \in S, \ u,v \in X^*$. Of course, by using
 Shirshov's Composition-Diamond Lemma, it can be easily seen that the set
 $Red(S)$ of $S-$reduced Lyndon-Shirshov words constitutes a linear
 basis of the quotient algebra $Lie(X)/ Id(S)$. The converse is
 also true.

  \ \ Explicitly Shirshov's Composition-Diamond Lemma for
  associative algebra was formulated by L. A. Bokut \cite{b76} in 1976 and
  G. Bergman \cite{b} in 1978.

   \ \ In this paper, we give  a comprehensive proof of Shirshov's Composition-Diamond
   Lemma for Lie algebras. There is an elementary approach to
   Gr\"{o}bner-Shirshov bases theory,
   including for Lie algebras, in \cite{bs}. We use properties of associative
   Lyndon-Shirshov words (ALSW) and non-associative
   Lyndon-Shirshov words (NLSW), see for example, \cite{L}. These properties are
   found by using
   induction on the length of a word applying Shirshov's elimination
   procedure of \cite{S53} (it is known also as the Lazard or Lazard--Shirshov
elimination, cf. \cite{R} and \cite{L}).

    \ \ This paper is based on the lectures given by the first author at
   Novosibirsk State University in September-October, 2006. Some
   notes were rewritten in a seminar at South China Normal University
   for master degree students. We thank Mr. Yu Li and Ms.
   Hongshan Shao for many valuable comments.

\section{Preliminaries}
We start with the Lyndon-Shirshov associative words.

 Let $X=\{x_i
|i\in I\}$ be a well-ordered set with $x_i>x_p$ if $i>p$ for any
$i,p\in I$. Let $X^*$ be the free monoid generated by $X$. For
$u=x_{i_1}x_{i_2} \cdots x_{i_k} \in X^*$, let
\begin{eqnarray*}
x_\beta=min(u)&=&min\{x_{i_1},x_{i_2}, \cdots ,x_{i_k}\}, \\
fir(u)&=&x_{i_1},\\
length\ of\ u:\ |u|&=&k.
\end{eqnarray*}

\begin{definition}
 Let $u=x_{i_1}x_{i_2} \cdots x_{i_k}\in X^*$. Then $u$ is called Weak-ALSW  if
 $fir(u)>min(u)$ or $|u|=1$, where ALSW means an ``associative Lyndon-Shirshov word".
\end{definition}
Let $u$ be a Weak-ALSW, $ min(u)=x_\beta$ and $|u|\geq2$. We define
$$
X'(u)=\{x_i^j=x_i\underbrace{x_\beta \cdots x_\beta}_{j}| i>\beta,
j\geq0\}.
$$
Note that $x_i^j=x_i\underbrace{x_\beta \cdots x_\beta}_{j}$ is just
a symbol.

 Now, we order $X'(u)$ by the following way:
$$
x_{i_1}^{j_1}>x_{i_2}^{j_2} \ \Leftrightarrow i_1>i_2 \ or \
(i_1=i_2, \ j_2>j_1).
$$
Suppose that $u,v$ are Weak-ALSW's and $min(v)\geq min(u)=x_\beta$.
Then we define
$$
v'_u=x_{i_1}^{m_1} \cdots x_{i_t}^{m_t} \ in \ (X'{(u)})^*
\Leftrightarrow v=x_{i_1}\underbrace{x_\beta \cdots x_\beta}_{m_1}
\cdots x_{i_t}\underbrace{x_\beta \cdots x_\beta}_{m_t} \ in \ X^*,
$$
where $x_{i_j}>x_\beta$, $m_j\in N$, $1\leq j\leq t$.  For the sake
of simpler
notation, we use $u'$ instead of $u'_{u}$. \\

Throughout Section 2 and 3, we assume that $x_1<x_2<x_3<\cdots$.

\begin{example}
Let $u=x_2x_1, v=x_3x_2$.  Then
 $v'_u=x^0_3x^0_2,  \
\ v'=x^1_3$.
\end{example}

The following lemma is obvious.
\begin{lemma}\label{l2.3}
Let u be a Weak-ALSW, $x_\beta=min(u), \ \  u=vw, \ v,w\neq1$ and
$w\neq
 x_\beta w_1 $. Then $u'=v'_uw'_u$.
\end{lemma}

\begin{example}
Let $v=x_3x_2x_1, w=x_2x_2 \  and \  u=vw=x_3x_2x_1x_2x_2.$ Then
$u'=x^0_3x^1_2x^0_2x^0_2$, $v'_u=x^0_3x^1_2,$ $w'_u=x^0_2x^0_2$ and
 $u'=v'_uw'_u.$
\end{example}

Recall that without specific explanation, we always use the
lexicographic order both on $(X'{(u)})^*$ and $X^*$ (i.e., $w>wt$ if
$t\neq 1$ and $zx_it_1>zx_jt_2$ if $x_i>x_j$).

\begin{lemma}\label{l2.5}
Let $u,v$ be Weak-ALSW's with $|v|\geq 2$. Then $u>v\Leftrightarrow
u'_{uv}>v'_{uv}.$
\end{lemma}
{\it Proof.} Let $x_\beta=min(uv)$. Assume that $u>v$. Then there
are two cases to consider.

 Case 1:\
\begin{eqnarray*}
u&=&x_{i_1}\underbrace {x_ \beta \cdots x_ \beta}_{l_1} \cdots
x_{i_{s-1}} \underbrace {x_ \beta \cdots x_
\beta}_{l_{s-1}}x_{i_{s}}\underbrace
{x_ \beta \cdots x_ \beta}_{l_s} \cdots\\
v&=&x_{i_1}\underbrace {x_ \beta \cdots x_ \beta}_{l_1} \cdots
x_{i_{s-1}} \underbrace {x_ \beta \cdots x_ \beta}_{l_{s-1}}yz
\cdots,\ \ \
 where \ \ x_{i_s}>y.
 \end{eqnarray*}
 \begin{enumerate}
 \item[(a)] If $y=x_ {\beta}$, then

\begin{eqnarray*}
u'_{uv}&=&x_{i_1}^{l_1} \cdots
x_{i_{s-2}}^{l_{s-2}}x_{i_{s-1}}^{l_{s-1}}x_{i_{s}}^{l_{s}}\cdots\\
v'_{uv}&=&x_{i_1}^{l_1} \cdots
x_{i_{s-2}}^{l_{s-2}}x_{i_{s-1}}^{l_{s-1}'}\cdots, \ \ where \ \
{l_{s-1}'}>l_{s_1}.
\end{eqnarray*}
So, $u'_{uv}>v'_{uv}$.

\item[(b)] If $y>x_ {\beta}$, then
\begin{eqnarray*}
u'_{uv}&=&x_{i_1}^{l_1} \cdots
x_{i_{s-1}}^{l_{s-1}}x_{i_{s}}^{l_{s}}\cdots\\
v'_{uv}&=&x_{i_1}^{l_1} \cdots x_{i_{s-1}}^{l_{s-1}} y^n \cdots. \ \
\end{eqnarray*}
So, $u'_{uv}>v'_{uv}$.
\end{enumerate}
Case 2:\
\begin{eqnarray*}
u&=&x_{i_1}\underbrace {x_ \beta \cdots x_ \beta}_{l_1} \cdots
x_{i_{s}}\underbrace
{x_ \beta \cdots x_ \beta}_{l_s} \\
v&=&x_{i_1}\underbrace {x_ \beta \cdots x_ \beta}_{l_1} \cdots
x_{i_{s}} \underbrace {x_ \beta \cdots x_ \beta}_{l_{s}}yz \cdots.
\end{eqnarray*}
\begin{enumerate}
 \item[(a)] If $y=x_ {\beta}$, then
\begin{eqnarray*}
u'_{uv}&=&x_{i_1}^{l_1} \cdots
x_{i_{s}}^{l_{s}}\\
v'_{uv}&=&x_{i_1}^{l_1} \cdots x_{i_{s}}^{l_{s}'}\cdots, \ \ where \
\ {l_{s}'}>l_{s}.
\end{eqnarray*}
So, $u'_{uv}>v'_{uv}$.
\item[(b)] If $y>x_ {\beta}$, then
$v'_{uv}=u'_{uv}y^n\cdots$ and so, $u'_{uv}>v'_{uv}$.
\end{enumerate}

Conversely,  assume that $u'_{uv}>v'_{uv}$. We will prove that
$u>v$. There are also two cases to consider. \\

Case 1:\ $u'_{uv}=x_{i_1}^{l_1}\cdots x_{i_s}^{l_s}, \ \
v'_{uv}=x_{i_1}^{l_1}\cdots x_{i_s}^{l_s}y^n\cdots.$ \\

Case 2:\ $u'_{uv}=x_{i_1}^{l_1}\cdots
x_{i_{s-1}}^{l_{s-1}}x_{i_s}^{l_s}\cdots, \ \
v'_{uv}=x_{i_1}^{l_1}\cdots
x_{i_{s-1}}^{l_{s-1}}{x_{i_s'}}^{{l_s}'}\cdots $, \  where
$x_{i_s}>x_{i_s'}$  or $(x_{i_s}=x_{i_s'}\ \ and \ \
l_s'>l_s)$. \\

In both cases, it is clear that $u>v$.\ \ \ \ $\square$

\begin{definition}
Let $u \in X^*$. Then $u$ is called  an ALSW if
$$
(\forall  v,w \in X^*, \ v,w\neq 1) \ \  u=vw\Rightarrow vw>wv.
$$
\end{definition}
\noindent{\bf Remark:} Let $u,v\in X^*$ and the $v_u'\in (X'(u))^*$
be as before. We denote by $|v|$ the length of $v$ in $X^*$ and
$|v_u'|_{X'}$ the length of $v_u'$ in $(X'(u))^*$.

\begin{lemma}\label{l2.7}
Let $u$ be a  Weak-ALSW with $|u|\geq2$. Then $u$ is an ALSW in
$X^*$ if and only if  $u'$  is an ALSW in $(X'{(u)})^*$.
\end{lemma}
{\it Proof.}  $``\Longrightarrow"$ \ \ If $|u'|_{X'}=1$, then $u'$
is an ALSW. Suppose that $|u'|_{X'}>1$ and $u'=v'_uw'_u$. Then
$u=vw$. Since $u$ is
 an ALSW, $vw>wv$ which implies  $(vw)'_u>(wv)'_u$ by Lemma 2.5.
 Therefore, by Lemma 2.3,
 $v'_u w'_u>w'_u v'_u$ and so, $u'$ is an ALSW.

 \ \ $``\Longleftarrow"$ \ \ Let $u=vw$ and $x_\beta=min(u)$. \ If \
$fir(w)=x_\beta$, \ then $vw>wv$. \ If \ $fir(w)\neq x_\beta$, then
$$
u'=v'_uw'_u\ \Rightarrow\ v'_uw'_u>w'_uv'_u\ \Rightarrow\
(vw)'_u>(wv)'_u\ \Rightarrow\ vw>wv.
$$
Hence, $u$ is an ALSW.

\noindent{\bf Remark}: For a Weak-ALSW  $u$, it is clear that
$|u'|_{X'}<|u|$ if $|u|>1$. For an ALSW $u$, we denote by $u''=
(u')'$ and $u^{(k)}=(u')^{(k-1)}$ for $k>0$ generally. From this, it
follows that $X^k(u)=X^{k-1}(u')$.

\begin{lemma}\label{l2.8}
For $u \in X^*$, $u$ is an ALSW if and only if $( \exists k \geq
0)$, s.t., \ $|u^{(k)}|_{X^{k}{(u)}}=1$.
\end{lemma}
{\it Proof.}
 We apply induction on $|u|$.  If $|u|=1$, then there is
 nothing to do.
Assume that $|u|>1$. Since $|u'|_{X'}<|u|$ and
$$
|u^{(k)}|_{X^k(u)}=|(u')^{(k-1)}|_{X^{k-1}(u')},
$$
by induction and by Lemma 2.7, the result follows.\ \ \ $\square$

\begin{example}
Let $u=x_5x_4x_5x_3$. Then
$$
u'=x^0_5x^0_4x^1_5, \ u''=(x^0_5)^1(x^1_5)^0 \ and \
u'''=((x^0_5)^1)^1.
$$
Therefore, by Lemma 2.8, $u$ is an ALSW.
\end{example}
\begin{lemma}\label{l2.10}
Let $u \in X^*$. Then $u$ is an ALSW if and only if
$$
(\forall  v,w \in X^*, \ v,w\neq 1) \ \  u=vw\Rightarrow u>w.
$$
\end{lemma}
{\it Proof.} \ \ $``\Longrightarrow"$ \ Induction on $|u|$. If
$|u|=1$, then the result clearly holds. Suppose that $ |u|\geq 2 , \
x_\beta=min(u)$ and $u=vw, \ v,w\neq 1$. If $w=x_ \beta w_1$, then
$u>w$. If $w\neq x_ \beta w_1$, then $u'=v'_uw'_u$. Since $u'$ is an
ALSW, by induction, $u'>w'_u$. Hence, by Lemma 2.5, $u>w$.

 \ \ $``\Longleftarrow"$ \ Induction on $|u|$. If $|u|=1$, then $u=x_i$ is an ALSW. If
$|u|>1$ and $|u'|_{X'}=1$, then by Lemma 2.8, $u$ is an ALSW. If
$|u'|_{X'}>1$ and $u'=v'_uw'_u$, then $u'>w'_u$ follows from $u>w$.
By induction, $u'$ is an ALSW. Hence, by Lemma 2.7, $u$ is an ALSW.\
\ \ \ $\square$

\begin{lemma}\label{l2.11}
Suppose that $u$ is an ALSW, $x_\beta=min(u)$ and $|u|>1$. Then
$ux_\beta$ is an ALSW.
\end{lemma}
{\it Proof.} Follows from Lemma 2.10.\ \ \ \ $\square$
\begin{lemma}\label{l2.12}
Let $u$ and $v$ be ALSW's. Then $uv$ is an ALSW  \ if and only if \
$u>v$.
\end{lemma}
{\it Proof.} \ $``\Longrightarrow"$ Suppose that $uv$ is an ALSW.
Then, by Lemma 2.10, $u>uv>v$.

 \ $``\Longleftarrow"$ \ We use induction on $|uv|$. Suppose that $u>v$. If $|uv|=2$
  or $v=x_{\beta}=min(uv)$, then
 the result is obvious. Otherwise, we can get that $u'_{uv}>v'_{uv}$,
 where $u'_{uv},v'_{uv}$ are ALSW's. By induction, $u'_{uv}v'_{uv}=(uv)'$ is an
 ALSW and so is $uv$. \ \ \ \ $\square$

\begin{lemma}\label{l2.13}
For any $ u \in X^*$, there exists a unique decomposition $u=u_1u_2
\cdots u_k$, where $u_i$ is an ALSW, \ $1 \leq i \leq k$, \ \ and \
$u_1 \leq u_2 \leq \cdots \leq u_k$.
\end{lemma}
{\it Proof.} \ To prove the existence, we use induction on $|u|$. If
$|u|=1$ then it is trivial. Let $|u|>1$ and $x_{\beta}=min(u)$. If
$u=x_{\beta} v$, then $v$ has the required decomposition and so does
$u$. Otherwise, $u$ is a Weak-ALSW. Thus, $u'$ has the decomposition
and so does $u$, by Lemma 2.5 and Lemma 2.7.

To prove the uniqueness, we let $ u=u_1\cdots u_k=w_1\cdots w_s $ be
the decompositions such that $u_i,w_j$ are ALSW's for any $i,j$;
$u_1\leq\cdots \leq u_k$ and $w_1\leq\cdots\leq w_s$. If
$u=x_{\beta} v$, then $u_1=w_1=x_{\beta}$ and the result follows
from the induction on $|u|$. Otherwise, $u$ is a Weak-ALSW and
$u'=u_{1u}'\cdots u_{ku}'=w_{1u}'\cdots w_{su}'$ are the
decompositions of $u'$. Now, by induction again, the result follows.
\ \ \ \ $\square$

\noindent{\bf Remark:} In Lemma 2.13, the word $u_k$ is the longest
ALSW end of $u$.

\begin{example} Let $u=x_1x_1x_2x_1x_2x_1x_1$. Then
$$
u=\underbrace{x_1}_{u_1}\underbrace{x_1}_{u_2}\underbrace{x_2x_1x_2x_1x_1}_{u_3}
=u_1u_2u_3
$$
is the decomposition of $u$. \ \
\end{example}

\begin{lemma}\label{l2.15}
Let $u$ be an ALSW \ and $|u| \geq 2$. If $u=vw$, where $w$ is the
longest ALSW proper end of $u$, then $v$ is an ALSW.
\end{lemma}
{\it Proof.} Suppose that $v$ is not an ALSW. Then, by Lemma 2.13,
we can assume that
$$
v=v_1v_2 \cdots v_m \ (m>1),
$$
where each $v_i$ is an ALSW and $v_1 \leq v_2 \leq \cdots \leq v_m$.
\ \ If $v_m>w$, then $v_mw$
 is an ALSW and $|v_mw|>|w|$, a contradiction.
If $v_m\leq w$, then we get another decomposition of $u$ which
contradicts the uniqueness in Lemma 2.13.  Thus, $v$ must be an
ALSW. \ \ \ \ $\square$

\begin{example} Let $u=x_5x_4x_5x_4x_3x_5x_3$. Then
$$
u=\underbrace{x_5x_4}_v\underbrace{x_5x_4x_5x_3}_w=vw
$$
and $u,v,w$ are all ALSW's.
\end{example}

Now, for an ALSW $u$, we introduce two bracketing ways.

One is up-to-down bracketing which is defined inductively by
$$
[x_i]=x_i, \ [u]=[[v][w]],
$$
where $u=vw$ and $w$ is the longest ALSW  proper end of $u$.

\begin{example} Let $u=x_2x_2x_1x_1x_2x_1$. Then
\begin{eqnarray*}
u&\rightarrow& [[x_2x_2x_1x_1][x_2x_1]]
\rightarrow[[x_2[x_2x_1x_1]][x_2x_1]]
 \rightarrow[[x_2[[x_2x_1]x_1]][x_2x_1]].
\end{eqnarray*}
\end{example}

The other is down-to-up bracketing. Let us explain it on a sample
word
$$
u=x_2x_2x_1x_1x_2x_1.
$$
Join the minimal letter $x_1$ to the previous letters:
$$
u\mapsto x_2[x_2x_1]x_1[x_2x_1].
$$
Form a new alphabet of the nonassociative words $x_2, \ [x_2x_1]$
and $x_1 $ ordered lexicographically, i.e.,
$$x_2>[x_2x_1]>x_1.
$$
Join the minimal letter $x_1$ to the previous letters:
$$
x_2[x_2x_1]x_1[x_2x_1] \mapsto  x_2[[x_2x_1]x_1][x_2x_1].
$$
Form a new alphabet
$$
x_2>[x_2x_1]>[[x_2x_1]x_1].
$$
Join the minimal letter $[[x_2x_1]x_1]$ to the previous letter:
$$
x_2[[x_2x_1]x_1][x_2x_1] \mapsto [x_2[[x_2x_1]x_1]][x_2x_1].
$$
Form a new alphabet
$$
[x_2[[x_2x_1]x_1]]>[x_2x_1].
$$
Finally, join the minimal letter $[x_2x_1]$ to the previous letter:
$$
[x_2[[x_2x_1]x_1]][x_2x_1] \mapsto [[x_2[[x_2x_1]x_1]][x_2x_1]]=[u].
$$

\noindent{\bf Remark:}  We denote by $[\ ]$ the down-to-up
bracketing and by $[[\ ]]$ the up-to-down bracketing.

\begin{lemma}
 $[u]=[[u]]$ for any ALSW $u$.
\end{lemma}
{\it Proof.} \ \ We use induction on $|u|$. If $|u|=1$, then $u=x_i$
and $[x_i]=[[x_i]]=x_i$. Assume that $u=vw$, where $w$ is the
longest ALSW proper end of $u$. Then, by Lemma 2.15, $v$ is an ALSW.
If $w=x_{\beta}w_1$, then
$$ w=x_ \beta \ and \ v=x_ix_\beta\cdots x_\beta.
$$
 By induction, we can get $[v]=[[v]]$. Hence, in this case, $$
[u]=[[u]]=([v]x_\beta).
$$
If $w\neq x_ \beta w_1$, then $u'=v'_uw'_u$. Suppose that
$v'_u={v_1}_{_{u}}'{v_2}_{{_u}}'$ such that ${v_2}_{{_u}}'w'_u$ is
an ALSW. Then $v_2w$ is an ALSW and $|v_2w|>|w|$, a contradiction.
This proves that $w'_u$ is the longest ALSW proper end of $u'$. By
induction, $[v'_u]=[[v'_u]]$ and $[w'_u]=[[w'_u]]$. Moreover, by the
definition of $[ \ ]$ and $[[\ ]]$, we have
$$
[ \ ]: \ u \mapsto v'_uw'_u \mapsto [v'_u][w'_u] \mapsto \cdots
$$
$$
[[\ ]]: \ u \rightarrow v'_uw'_u\rightarrow [[v'_u]][[w'_u]]
\rightarrow \cdots.
$$
Therefore, $[u]=[[u]]$. \ \ \ \ $\square$

\section{Free Lie algebras}

Now we give the definition of a non-associative Lyndon-Shirshov
word.

\begin{definition}\label{d3.1}
Let $<$ be the order on $X^*$ as before and $(u)$ a non-associative
word. Then $(u)$ is called a non-associative Lyndon-Shirshov word,
denoted by NLSW, if
\begin{enumerate}
\item[(i)] $u$ is an ALSW,
\item[(ii)] if $(u)=((v)(w))$, then both $(v)$ and $(w)$ are NLSW's,
\item[(iii)] in (ii) if $(v)=((v_1)(v_2))$, then $v_2 \leq w$ in $X^*$.
\end{enumerate}
\end{definition}

\noindent{\bf Remark}: In Definition 3.1 (ii), $v>w$ by Lemma 2.12.

\begin{theorem}
Let $u$ be an ALSW. Then there exists a unique bracketing way such
that $(u)$ is a NLSW.
\end{theorem}
{\it Proof.} (Existence). \ Let $u$ be an ALSW. We will prove that
up-to-down bracketing is one of bracketing way such that $[[u]]$ is
a NLSW. Induction on $|u|$. If $|u|=1$, then nothing to do.
  Suppose
that $|u|>1$ and $u=vw$ where $w$ is the longest ALSW proper end of
$u$. Then, $[[u]]=[[[v]][[w]]]$. By induction, both $[[v]]\ and \
[[w]]$ are NLSW's. Now, we assume that
$$
[[v]]=[[v_1]][[v_2]] \ and \ v_2>w.
$$
Then, $v_2 w$ is an ALSW, a contradiction. So, $v_2\leq w$ and
hence, $[[u]]$ is a NLSW.

(Uniqueness).\ We assume that $u$ is an ALSW and $( \ )$ is a
bracketing way such that $(u)$ is a NLSW. Then, we have to show $
(u)= [[u]] $. We use induction on $|u|$. If $|u|=1$, then
$(u)=[[u]]$ clearly. Suppose that
$$
|u|>1 \ \mbox{ and } \ u=x_{i_1}x_ \beta \cdots x_ \beta \cdots
x_{i_s} x_ \beta \cdots x_ \beta,
$$
where $x_{i_j}> x_ {\beta} = min(u)$.

Note that if $v=x_{i}x_ \beta \cdots x_ \beta,\ x_{i}>x_ \beta $,
then
$$
[[v]]=[[\cdots[[x_{i}x_ \beta]] \cdots x_\beta]] x_ \beta
$$
is the unique bracketing way such that $[[v]]$ is a NLSW. According
to the definition of NLSW, any associative word in a bracket must be
ALSW. Hence,
\begin{eqnarray*}
(u)&=&((x_{i_1}x_ \beta \cdots x_ \beta)(x_{i_2}x_ \beta \cdots x_
\beta)\cdots (x_{i_s}x_ \beta \cdots x_ \beta)) \\
&=& [[[[x_{i_1}x_ \beta \cdots x_ \beta]][[x_{i_2}x_ \beta \cdots x_
\beta]]\cdots [[x_{i_s}x_ \beta \cdots x_ \beta]]]].
\end{eqnarray*}
By induction, $(u')=[[u']]$ and therefore, $[[u]]=[[u']]=(u')=(u)$.
\ \ \ \ $\square$

 Let $X^{**}$ be the set
of all non-associative words $(u)$ in $X$. If $(u)$ is a NLSW, then
we denote it by $[u]$.

From now on, let $k\langle X\rangle$ be the free associative algebra
generated by $X$. We consider $(\ )$ as Lie bracket in $k\langle
X\rangle$, i.e., for any $a,b\in k\langle X\rangle, \ (ab)=ab-ba$.
Denote by $Lie(X)$ the subLie-algebra of $k\langle X\rangle$
generated by $X$.

Given a polynomial $f\in k\langle X\rangle$, it has the leading word
$\bar f \in X^*$ according to the above order on $X^*$ such that
$$
f=\sum_{u\in X^*}f(u)u=\alpha\overline{f}+\sum{\alpha}_iu_i,
$$
where $\overline{f}, \ u_i\in X^*$, $\overline{f}>u_i, \ \alpha , \
{\alpha}_i, \ f(u)\in k$. We call $\overline{f}$ the leading term of
$f$. Denote the set $\{ u|f(u)\neq 0\}$ by $suppf$ and $deg(f)$ by
$|\overline{f}|$. $f$ is called monic if $\alpha=1$.

Note that if $|u|=|v|$ and $u<v$, then the lexicographic order which
we use before is the same as the degree-lexicographic order on
$X^*$.

\begin{theorem}\label{t3.3}
Let the order $<$ be as before. Then, for any $(u) \in X^{**}, \
(u)$ has a representation:
$$
(u)=\sum \alpha_i [u_i],
$$
where each $\alpha_i \in k, \ [u_i]$ is a NLSW  and $|u_i|=|u|$.
Even more, if $(u)=([v][w])$, then $u_i>min\{v,w\}$.
\end{theorem}
{\it Proof.} Induction on $|u|$. If $|u|=1$, then $(u)=[u]$ and the
result holds. Suppose that $|u|>1$ and $(u)=((v)(w))$. Then, by
induction,
$$
(v)=\sum \alpha_i [v_i] \ \mbox{ and }\ (w)=\sum \beta_j [w_j],
$$
where $\alpha_i,\beta_j \in k, \ [v_i], [w_j]$ are NLSW's,
$|v_i|=|v|$ and $|w_j|=|w|$. Without loss of generality, we may
assume that $(u)=([v][w])$ with $v>w$ because of
$([v][w])=-([w][v])$. If $|v|=1$, then
$$
(u)=([v][w])
$$
is a NLSW. Suppose that $|v|>1$ and $[v]=[[v_1][v_2]]$.

 There are
two subcases
\begin{enumerate}
\item[(a)] If $v_2\leq w$, then $(u)=(([v_1][v_2])[w])$ is a NLSW.
\item[(b)] If $v_2>w$, then
$$
(u)=(([v_1][v_2])[w])=(([v_1][w])[v_2])+([v_1]([v_2][w])).
$$
\end{enumerate}
By induction,
\begin{eqnarray*}
([v_1][w])&=&\sum \gamma_i [t_i],\ \ t_i>min\{v_1,w\}=w\\
([v_2][w])&=&\sum {\gamma_j}' [t_j'],\ \ {t_j'}>min\{v_2,w\}=w.
\end{eqnarray*}
Then,
$$
(u)=\sum \gamma_i ([t_i][v_2])+ \sum {\gamma_j}' ([v_1][t_j']).
$$
By noting that
$$
min\{t_i,v_2\} \ \mbox{ and } \ min\{t_j',v_1\}>min\{v,w\}=w,
$$
the result follows from the inverse induction on $min\{v,w\}$. \ \ \
\ $\square$

\begin{example} Let $(u)=(((x_3x_2)(x_2x_1))(x_2x_1x_1))$. Then
\begin{eqnarray*}
(u)&=&(((x_3(x_2x_1))x_2)(x_2x_1x_1))+((x_3(x_2(x_2x_1)))(x_2x_1x_1)),\\
(((x_3(x_2x_1))x_2)(x_2x_1x_1))&=&(((x_3(x_2x_1))(x_2x_1x_1))x_2)+((x_3(x_2x_1))(x_2(x_2x_1x_1)))\\
&=&(((x_3(x_2x_1x_1))(x_2x_1))x_2)+((x_3((x_2x_1)(x_2x_1x_1)))x_2)\\
& &+((x_3(x_2x_1))(x_2(x_2x_1x_1))),\\
((x_3(x_2(x_2x_1)))(x_2x_1x_1))&=&((x_3(x_2x_1x_1))(x_2(x_2x_1)))+(x_3((x_2(x_2x_1))(x_2x_1x_1)))\\
 &=&((x_3(x_2x_1x_1))(x_2(x_2x_1)))+(x_3((x_2(x_2x_1x_1))(x_2x_1))\\
& &+(x_3(x_2((x_2x_1)(x_2x_1x_1)))),
\end{eqnarray*}
and hence,
\begin{eqnarray*}
(u)&=&(((x_3(x_2x_1x_1))(x_2x_1))x_2)+((x_3((x_2x_1)(x_2x_1x_1)))x_2)\\
& &+((x_3(x_2x_1))(x_2(x_2x_1x_1)))+((x_3(x_2x_1x_1))(x_2(x_2x_1)))\\
& &+(x_3((x_2(x_2x_1x_1))(x_2x_1))+(x_3(x_2((x_2x_1)(x_2x_1x_1))))
\end{eqnarray*}
is a linear combination of NLSW's.
\end{example}

\begin{lemma}\label{l3.5}
Let $[u]$ be a NLSW. Then $\overline{[u]}=u$.
\end{lemma}

{\it Proof.} We use induction on $|u|$. If $|u|=1$, then the result
holds immediately. Let $|u|>1$ and $[u]=[[v][w]]$. Then, by
induction, $\overline{[v]}=v$ and $\overline{[w]}=w$. Suppose that
$$
[v]=v+ \sum \limits_{v_i<v} \alpha_iv_i, \ \ [w]=w+ \sum
\limits_{w_j<w} \beta _jw_j,
$$
where $\alpha_i,\beta_j\in k, \ v,v_i,w,w_j\in X^*$. It is easy to
see that $|v_i|=|v|$ and $|w_j|=|w|$ for any $i,j$. Then,
\begin{eqnarray*}
[u]&=&[(v+ \sum \limits_{v_i<v} \alpha _iv_i)(w+ \sum
\limits_{w_j<w} \beta _jw_j)] \\
&=&(v+ \sum \limits_{v_i<v} \alpha _iv_i)(w+ \sum \limits_{w_j<w}
\beta _jw_j)-(w+ \sum \limits_{w_j<w} \beta
_jw_j)(v+ \sum \limits_{v_i<v} \alpha _iv_i)\\
&=&vw + \sum \limits_{w_j<w} \beta _jvw_j + \sum \limits_{v_i<v}
\alpha _iv_iw +\sum \alpha _i \beta_j v_i w_j\\
&-& wv - \sum \limits_{w_j<w}\beta_jw_jv - \sum \limits_{v_i<v}
\alpha _iwv_i - \sum \beta_j \alpha_i w_jv_i.
\end{eqnarray*}
Since
$$
vw>vw_j, \ v_iw, \ v_iw_j, \  wv \ \mbox{ and } \ wv>wv_i, \ w_jv,\
w_jv_i,
$$
we have, $\overline{[u]}=u$. \ \ \ \ $\square$

\noindent{\bf Remark.} By the proof of Lemma 3.5, if we consider
$[u]$ as a polynomial in $k\langle X\rangle$, then each $r\in
supp([u])$ has the same length as $u$, moreover, $cont(r)=cont(u)$,
where, for example, $cont(u)=\{x_{i_1},\cdots, x_{i_t}\}$ if
$u=x_{i_1}\cdots x_{i_t}\in X^*$.

\begin{lemma}\label{l3.6}
NLSW's are $k-$independent.
\end{lemma}
{\it Proof.}  Suppose
$$
\sum \limits_{i=1}^k \alpha _i[u_i]=0,
$$
where each $\alpha_i \in k$, $[u_i]$ is a NLSW and $u_1>u_2> \cdots
> u_k$. If $\alpha_1\neq 0$, then
$\overline{\sum\limits_{i}\alpha_i[u_i]}=u_1\neq 0$,  a
contradiction. Then, all $\alpha_i$ must be 0. \ \ \ \ $\square$

\ \

By Theorem 3.3 and Lemma 3.6, we have the following corollary.

\begin{corollary}\label{c3.7}
NLSW's are linear basis of  $Lie(X)$.
\end{corollary}

From Corollary 3.7 and Lemma 3.5, we have

\begin{corollary}
For any $ f \in Lie(X)$, $\overline{f}$ is an ALSW.
\end{corollary}

\begin{theorem}
$Lie(X)$ is the free Lie algebra generated by $X$.
\end{theorem}
{\it Proof.} Let $L$ be a Lie algebra and $f:\ X\longrightarrow \ L$
a mapping. Then, we define a mapping
$$
\bar f :\ Lie(X)\longrightarrow L;\ [x_{i_1}\cdots
x_{i_n}]\longmapsto [f(x_{i_1})\cdots f(x_{i_n})],
$$
where $[x_{i_1}\cdots x_{i_n}]$ is NLSW. It is easy to check $\bar
f$ is a unique Lie homomorphism such that $\bar f i=f$.

\setlength {\unitlength}{1cm}
\begin{picture}(7, 3)
\put(4.2,2.3){\vector(1,0){1.7}} \put(4.1, 2.0){\vector(0,-1){1.3}}
\put(5.9,2.1){\vector(-1,-1){1.6}} \put(3.9,0.2){$L$}
\put(6,2.2){$Lie(X)$} \put(3.9,2.2){$X$} \put(4.9, 2.4){$i$}
\put(3.7,1.3){$f$} \put(5.3,1){$\exists!\ \overline{f}$}
\end{picture}\ \ \ \ \ $\square$

 The
following theorem plays a key role in proving the
Composition-Diamond lemma for Lie algebras (see Theorem 5.8).

\begin{theorem}\label{t4} {\em (A. I. Shirshov \cite{S58})}\ Let $u,v$ be ALSW's,
$u=avb, \ a,b\in{X^*}$. Then
\begin{enumerate}
\item[(i)] $[u]=[a[vc]d]$, where $b=cd, \ c,d\in{X^*}$.\\
\item[(ii)] Let
\begin{equation}\label{e1}
[u]_v=[u]|_{[vc]\mapsto{[[[v][c_1]]\cdots[c_k]]}},
\end{equation}
where $c=c_1 \cdots c_k$, $c_j$ is an ALSW and $c_1 \leq{c_2
}\leq\cdots \leq{c_k}$. Then, in $k\langle X\rangle$,
$$
\overline{[u]_v}=u.
$$
Moreover, $$ [u]_v=a[v]b+\sum\limits_{i}\alpha_ia_i[v]b_i,$$ where
each $\alpha_i\in{k}$ and $a_ivb_i<avb$.
\end{enumerate}
\end{theorem}

{\it Proof.} $(i)$ Induction on $|u|$. If $|u|=1$, then $u=v=x_i$
and the result holds. Assume that $|u|>1$. If $v=x_i$, then
$[u]=[a[x_i]d]$ and the result holds.
 Now, we
consider the case of $|v|>1$. Let $x_\beta=min(u)$ and
$b=x_\beta^e\tilde{b}$, where $e\geq0$ and
$fir(\tilde{b})\neq{x_\beta}$. Then
$$
u=avb=avx_\beta^e\tilde{b}=a\tilde{v}\tilde{b},
$$
where $\tilde{v}=vx_\beta^e$ is also an ALSW, by Lemma 2.11. Then,
by induction, for $u'=a_u'\tilde{v}_u'\tilde{b}_u'$, we have
$[u']=[a_u'[\tilde{v}_u'\tilde{c}_u']d_u'], \
\tilde{b}_u'=\tilde{c}_u'd_u'$. By substitution \\
$x_i^j\mapsto[[x_ix_\beta]\cdots x_\beta]$, we obtain
$$
[u]=[a[\tilde{v}\tilde{c}]d]=[a[vx_\beta^e\tilde{c}]d]=[a[vc]d], \
\mbox{ where } c=x_\beta^e\tilde{c}.$$

$(ii)$ If $c=1$, then $[u]_v=[u]$ and the results hold clearly.
Otherwise, by Lemma 2.13, we may assume that
$$
c=x_\beta \cdots x_\beta c_{l+1}\cdots c_k,
$$
where each $c_i$ is an ALSW and $x_\beta<c_{l+1}\leq \cdots \leq{c_k}$. \\
Then
$$
[u]_v=[u]|_{[vx_\beta^e\tilde{c}]\mapsto{[[[v]x_\beta]\cdots
x_\beta[c_{l+1}]\cdots[c_k]]}} \ \mbox{ and }
[u]_{\tilde{v}}=[u]|_{[\tilde{v}\tilde{c}]\mapsto{[[[\tilde{v}][c_{l+1}]]\cdots[c_k]]}}.
$$
Now, we use induction on $|u|$. If $|u|=1$, then this is a trivial
case. Suppose that $|u|>1$ and $|v|>1$. Then, by (i),
$$
u=a\tilde{v}\tilde{c}d, \ \ \ u'=a_u'\tilde{v}_u'\tilde{c}_u'd_u'
$$
and by induction,
$$
[u']_{\tilde{v}_u'}=a_u'[\tilde{v}_u']\tilde{c}_u'd_u'+\sum\limits_{i\in{I_1}}\alpha_i{a_i}
_u'[\tilde{v}_u']{b_i} _u',
$$
where each ${a_i}_u'\tilde{v}_u'{b_i}_u'<u'$. Now, it is easy to
check that
$$
[[x_ix_\beta]\cdots x_\beta]=\sum \limits_{m\geq0}(-1)^m\left(
\begin{array}{c}
j \\
m \\
\end{array}
\right)x_\beta^mx_ix_\beta^{j-m}\ \mbox{ and }\
x_\beta^mx_ix_\beta^{j-m}<x_ix_\beta^j\ \ (m>0).
$$
Now, by substitution $ x_i^j\mapsto[[x_ix_\beta]\cdots x_\beta]$, we
obtain
$$
[u]_{\tilde{v}}=a[\tilde{v}]\tilde{c}d+\sum\limits_{i\in{I_2}}\alpha_ia_i[\tilde{v}]b_i,
$$
where each $a_i\tilde{v}b_i<a\tilde{v}\tilde{c}d$. Also, by
substitution $ [\tilde{v}]\mapsto{[[[v]x_\beta]\cdots x_\beta]}, $
we have
$$
[u]_v=a[v]x_\beta^e\tilde{c}d+\sum\limits_{j\in{I}}\beta_ja_j[v]b_j
=a[v]b+\sum\limits_{j\in{I}}\beta_ja_j[v]b_j,
$$
where each $a_jvb_j<avb$. \ \ \ \ $\square$

\noindent{\bf Remark.} By the proof of Theorem 3.10, if we consider
$[u]_v$ as a polynomial in $k\langle X\rangle$, then for any $w\in
supp([u]_v)$, $cont(w)=cont(u)$.

\begin{definition}
Let $S\subset Lie(X)$ with each $s\in S$ monic, $a,b\in{X^*}$ and
$s\in S$. If $a\bar{s}b$ is an ALSW, then we call
$[asb]_{\bar{s}}=[a\bar{s}b]_{\bar{s}}|_{[\bar{s}]\mapsto{s}}$ a
normal $S$-word (or normal $s$-word) while $[a\bar{s}b]_{\bar{s}}$
is called a relative nonassociative Lyndon-Shirshov word, denoted by
RNLSW, where $[a\bar{s}b]_{\bar{s}}$ is defined by (3.1) (see
Theorem 3.10).
\end{definition}
\begin{corollary}\label{c3.11} Let $u,v$ be ALSW's, $f\in{Lie(X)}$,
$\bar{f}=v$ and $u=avb, \ a,b\in{X^*}.$  Then, for the normal
$f$-word $[afb]_v=[avb]_v|_{[v]\mapsto{f}}$, we have
$$
[afb]_v=afb+\sum\limits_{i}\alpha_ia_ifb_i,
$$
where each $\alpha_i\in{k}, \ a_i, b_i\in{X^*},\ a_i\bar{f}b_i<u$.\
\end{corollary}

\section{Composition-Diamond lemma for associative algebras}
In this Section, we cite some concepts and results from the
literature which are related to the Gr\"{o}bner-Shirshov basis for
the associative algebras.
\begin{definition}  {\em (\cite{Sh}, see also \cite{b72}, \cite{b76})} \
Let $f$ and $g$ be two monic polynomials in \textmd{k}$\langle
X\rangle$ and $<$ a well order on $X^*$. Then, there are two kinds
of compositions:
\begin{enumerate}
\item[(i)] If  $w$ is a word such that $w=\bar{f}b=a\bar{g}$ for some
$a,b\in X^*$ with deg$(\bar{f})+$deg$(\bar{g})>$deg$(w)$, then the
polynomial  $(f,g)_w=fb-ag$ is called the intersection composition
of $f$ and $g$ with respect to $w$.
\item[(ii)] If  $w=\bar{f}=a\bar{g}b$ for some $a,b\in X^*$, then the
polynomial $(f,g)_w=f - agb$ is called the inclusion composition of
$f$ and $g$ with respect to $w$.
\end{enumerate}
\end{definition}

\begin{definition} {\em (\cite{b72}, \cite{b76}, cf. \cite{Sh})}
Let $S\subset k\langle X\rangle$ with each $s\in S$ monic. Then the
composition $(f,g)_w$ is called trivial modulo $(S,w)$ if
$(f,g)_w=\sum\alpha_i a_i s_i b_i$, where each $\alpha_i\in k$,
$a_i,b_i\in X^{*}$ and $\overline{a_i s_i b_i}<w$. If this is the
case, then we write
$$
(f,g)_w\equiv_{ass}0\quad mod(S,w)
$$
In general, for $p,q\in k\langle X\rangle$, we write
$$
p\equiv_{ass} q\quad mod(S,w)
$$
which means that $p-q=\sum\alpha_i a_i s_i b_i $, where $\alpha_i\in
k,a_i,b_i\in X^{*}$ and $\overline{a_i s_i b_i}<w$.
\end{definition}

\begin{definition} {\em (\cite{b72}, \cite{b76}, cf. \cite{Sh})} \
We call the set $S$ with respect to the well order $<$ a
Gr\"{o}bner-Shirshov set (basis) in $k\langle X\rangle$ if any
composition of polynomials in $S$ is trivial modulo $S$.
\end{definition}

If a subset $S$ of $k\langle X\rangle$ is not a Gr\"{o}bner-Shirshov
basis, then we can add to $S$ all nontrivial compositions of
polynomials of $S$, and by continuing this process (maybe
infinitely) many times, we eventually obtain a Gr\"{o}bner-Shirshov
basis $S^{comp}$. Such a process is called the Shirshov algorithm.

A well order $>$ on $X^*$ is monomial if it is compatible with the
multiplication of words, that is, for $u, v\in X^*$, we have
$$
u > v \Rightarrow w_{1}uw_{2} > w_{1}vw_{2},  \mbox{ for  all }
 w_{1}, \ w_{2}\in  X^*.
$$
A standard example of monomial order on $X^*$ is the deg-lex order
to compare two words first by degree and then lexicographically,
where $X$ is a linearly ordered set.

The following lemma was proved by Shirshov \cite{Sh} for free Lie
algebras (with deg-lex ordering) in 1962 (see also Bokut
\cite{b72}). In 1976, Bokut \cite{b76} specialized the approach of
Shirshov to associative algebras (see also Bergman \cite{b}). For
commutative polynomials, this lemma is known as the Buchberger's
Theorem in \cite{bu65} and \cite{bu70}.

\begin{lemma}\label{l1}
{\em (Composition-Diamond Lemma)} \ Let $k$ be a field, $A=k \langle
X|S\rangle=k\langle X\rangle/Id(S)$ and $<$ a monomial order on
$X^*$, where $Id(S)$ is the ideal of $k \langle X\rangle$ generated
by $S$. Then the following statements are equivalent:
\begin{enumerate}
\item[(i)] $S $ is a Gr\"{o}bner-Shirshov basis.
\item[(ii)] $f\in Id(S)\Rightarrow \bar{f}=a\bar{s}b$
for some $s\in S$ and $a,b\in  X^*$.
\item[(ii')]
$f\in{Id(S)}\Rightarrow f=\alpha_1 a_1s_1b_1 +\alpha_2 a_2s_2b_2
+\cdots$, where $\alpha_i\in{k}$ and
$\bar{f}={a_1\bar{s_1}b_1}>{a_2\bar{s_2}b_2}>\cdots$.
\item[(iii)] $Red(S) = \{ u \in X^* |  u \neq a\bar{s}b ,s\in S,a ,b \in X^*\}$
is a basis of the algebra $A=k\langle X | S \rangle$.
\end{enumerate}
\end{lemma}

\section{Composition-Diamond lemma for Lie algebras}

 In this
section, we give the Composition-Diamond lemma for Lie algebras.

Throughout this section, we extend the lexicographic order on $X^*$
mentioned in Section 2 to the deg-lex order $<$ on $X^*$.

\begin{lemma}
Let $ac,cb$ be ALSW's, where $a,b,c\in X^*$ and $c\neq 1$. Then
$w=acb$ is also an ALSW.
\end{lemma}

{\it Proof.} We use induction on $|w|=n$. If $n=3$, then
$w=x_ix_jx_k$ is an ALSW, because $ac$ and $cb$ are ALSW's implies
that $x_i>x_j>x_k$. In the inductive case $n>3$, suppose
$min(w)=x_\beta$, $b=x_\beta^e\tilde{b},e\geq0, \ fir(\tilde{b})\neq
x_\beta$ and $\tilde{c}=cx_\beta^e$. Then
$$
w=a\tilde{c}\tilde{b} \ \mbox{ and }\
w'=a_w'\tilde{c}_w'\tilde{b}_w'.
$$
It is clear that $a_w'\tilde{c}_w',\tilde{c}_w'\tilde{b}_w'$ are
ALSW's. By induction, $w'$ is an ALSW and so is $w$. \ \ \ \
$\square$

\begin{definition}\label{d5.1}
Let $f$ and $g$ be two monic Lie polynomials in $Lie(X)\subset
k\langle X\rangle$. Then, there are two kinds of Lie compositions:
\begin{enumerate}
\item[(i)] If  $w=\bar{f}=a\bar{g}b$ for some $a,b\in X^*$, then the
polynomial $\langle f,g\rangle_w=f - [agb]_{\bar{g}}$ is called the
 composition of inclusion of $f$ and $g$ with respect to $w$.

\item[(ii)] If $w$ is a word such that $w=\bar{f}b=a\bar{g}$ for
some $a,b\in X^*$ with $deg$$(\bar{f})+$deg$(\bar{g})>$deg$(w)$,
then the polynomial
 $\langle f,g\rangle_w=[fb]_{\bar{f}}-[ag]_{\bar{g}}$ is called the composition of intersection of $f$ and
$g$ with respect to $w$.
\end{enumerate}
\end{definition}

By Lemma 5.1, in the Definition 5.2 (i) and (ii), $w$ is an ALSW.

\begin{definition}\label{d5.1}
Let $S\subset{Lie(X)}$ be a nonempty subset, $h$ a Lie polynomial
and $w\in X^*$. We shall say that $h$ is trivial modulo $(S,w)$,
denoted by $h\equiv_{Lie}0 \  mod(S,w)$, if
$h=\sum\limits_{i}\alpha_i[a_is_ib_i]_{\bar{s_i}}$, where each
$\alpha_{i}\in{k}$, $a_i, b_i\in{X^*}$, $s_i\in{S}$,
$[a_i\bar{s_i}b_i]_{\bar{s_i}}$ is a RNLSW and $a_i\bar{s_i}b_i<w$.
\end{definition}
\begin{definition} Let $S\subset{Lie(X)}$ be a nonempty set of
monic Lie polynomials. Then $S$ is called a Gr\"{o}bner-Shirshov set
(basis) in $Lie(X)$ if any composition $\langle f,g\rangle_w$ with
$f,g\in{S}$ is trivial modulo $(S,w)$, i.e., $\langle
f,g\rangle_w\equiv_{Lie}0 \ mod(S,w)$.
\end{definition}
\begin{lemma}\label{l5.5}
 Let $f,g$ be monic Lie polynomials. Then
$$
\langle f,g\rangle_w-(f,g)_w\equiv_{ass}0 \ \ mod(\{f,g\},w).
$$
\end{lemma}
{\it Proof.} If $\langle f,g\rangle_w$ and $(f,g)_w$ are
compositions of intersection, where $w=\bar{f}b=a\bar{g}$, then, by
Corollary 3.12, we may assume that
$$
\langle f,g\rangle_w=[fb]_{\bar{f}}-[ag]_{\bar{g}}
=fb+\sum\limits_{I_1}\alpha_ia_ifb_i-ag-\sum\limits_{I_2}\beta_ja_jgb_j,
$$
where $a_i\bar{f}b_i,a_j\bar{g}b_j < \bar{f}b=a\bar{g}=w$. It
follows that
$$
\langle f,g\rangle_w-(f,g)_w\equiv_{ass}0 \ \ mod(\{f,g\},w).
$$
Similarly, for the case of the compositions of inclusion, we have
the same conclusion. \ \ \ \ $\square$
\begin{theorem} {\em (\cite{b99},\cite{b99-})} \
Let $S\subset{Lie(X)}\subset{k\langle X\rangle}$ be a nonempty set
of monic Lie polynomials. Then $S$ is a Gr\"{o}bner-Shirshov basis
in $Lie(X)$  if and only if $S$ is a Gr\"{o}bner-Shirshov basis in
$k\langle X\rangle$.
\end{theorem}
{\it Proof.} Note that, by the definitions, for any $f,g\in S$, they
have composition in $Lie(X)$  if and only if so do in $k\langle
X\rangle$.

Suppose that $S$ is a Gr\"{o}bner-Shirshov basis in $Lie(X)$. Then,
for any composition $\langle f,g\rangle_w$, we have
$$
\langle
f,g\rangle_w=\sum\limits_{I_1}\alpha_i[a_is_ib_i]_{\bar{s_i}},
$$
where $[a_i\bar{s_i}b_i]_{\bar{s_i}}$ are RNLSW's and
$a_i\bar{s_i}b_i<w$. By Corollary 3.12,
$$
\langle f,g\rangle_w=\sum\limits_{I_2}\beta_jc_js_jd_j,
$$
where each $c_j\bar{s_j}d_j<w$. Thus, by Lemma 5.5, we get
$$
(f,g)_w\equiv_{ass}0 \ mod(S,w).
$$
Hence, $S$ is a Gr\"{o}bner-Shirshov basis in $k\langle X\rangle$.

Conversely, assume that $S$ is a Gr\"{o}bner-Shirshov basis in
$k\langle X\rangle$. Then, for any composition $\langle
f,g\rangle_w$ in $S$, by Lemma 5.5,  we obtain
$$
\langle f,g\rangle_w\equiv_{ass}(f,g)_w\equiv_{ass}0 \ \ mod(S,w).
$$
Therefore, we can assume, by Lemma 4.4, that
$$
\langle f,g\rangle_w=\sum\limits_{I_1}\alpha_ia_is_ib_i,
$$
where $a_i\bar{s_i}b_i<w$ and
$w>a_1\bar{s_1}b_1>a_2\bar{s_2}b_2>\ldots$. By noting that $\langle
f,g\rangle_w\in Lie(X), \ \overline{\langle
f,g\rangle_w}=a_1\bar{s_1}b_1$ is an ALSW which shows that
$[a_1\bar{s_1}b_1]_{\bar{s_1}}$ is a RNLSW. Let $h_1=\langle
f,g\rangle_w-\alpha_1[a_1s_1b_1]_{\bar{s_1}}$. Clearly,
$\overline{h_1}<\overline{\langle f,g\rangle_w}$. Then, by Corollary
3.12,  we have
$$
h_1\equiv_{ass}0 \ mod(S,w).
$$
Now, by induction on $\overline{\langle f,g\rangle_w}$, we have
$$
\langle
f,g\rangle_w=\sum\limits_{I_2}\alpha_i[c_is_id_i]_{\bar{s_i}},
$$
where each $[c_i\bar{s_i}d_i]_{\bar{s_i}}$ is a RNLSW  and
$c_i\bar{s_i}d_i<w$. This proves that $S$ is a Gr\"{o}bner-Shirshov
basis in $Lie(X)$. \ \ \ \ $\square$
\begin{lemma}\label{l5.7}
Let  $S\subset{Lie(X)}$ with each $s\in S$ monic. Let
$$
Red(S)=\{[u] \ | \ [u] \mbox{ is a } NLSW, \ u\neq{a\bar{s}b}, \
s\in{S},\ a,b\in{X^*}\}.
$$
Then, for any  $h\in{Lie(X)}$, $h$ has a representation:
$$h=\sum\limits_{{[u_i]\in{Red(S)}, \ u_i\leq\bar{h}}}\alpha_i[u_i]+
\sum\limits_{{s_j\in{S}, \
a_j\bar{s_j}b_j\leq\bar{h}}}\beta_j[a_js_jb_j]_{\bar{s_j}}.$$
\end{lemma}
{\it Proof.} We can assume that
$h=\sum\limits_{i}\alpha_{i}[u_{i}]$, where each $[u_i]$ is a NLSW,
$0\neq{\alpha_{i}\in{k}}$ and $u_{1}>u_{2}>\cdots$. If
$[u_1]\in{Red(S)}$, then let $h_{1}=h-\alpha_{1}[u_1]$. If
$[u_1]\not\in{Red(S)}$, then there exists $s\in{S}$ and
$a_1,b_1\in{X^*}$ such that $u_1=a_1\bar{s_1}b_1$. Now, let
$$
h_1=h-\alpha_1[a_1s_1 b_1]_{\bar{s_1}}\in{Lie(X)}.
$$
Hence, in both cases, we have $\bar{h_1}<\bar{h}$. Now, the result
follows from induction on $\bar{h}$. \ \ \ \ $\square$
\begin{theorem}\label{cdL}
Let $S\subset{Lie(X)}\subset{k\langle X\rangle}$ be nonempty set of
monic Lie polynomials. Let $Id_{Lie}(S)$ be the Lie-ideal of
$Lie(X)$ generated by $S$. Then the following statements are
equivalent.
\begin{enumerate}
\item[(i)] $S$ is a Gr\"{o}bner-Shirshov basis in
$Lie(X)$.
\item[(ii)] $f\in{Id_{Lie}(S)}\Longrightarrow{\bar{f}=a\bar{s}b}$, for
some $s\in{S}$ and $a,b\in{X^*}$.
\item[(ii')]
$f\in{Id_{Lie}(S)}\Longrightarrow
f=\alpha_1[a_1s_1b_1]_{\bar{s_1}}+\alpha_2[a_2s_2b_2]_{\bar{s_2}}+\cdots$,
where $\alpha_i\in{k}$ and
$\bar{f}={a_1\bar{s_1}b_1}>{a_2\bar{s_2}b_2}>\cdots$.
\item[(iii)]$Red(S)=\{[u] \ | \ [u] \mbox{ is a } NLSW, \ u\neq{a\bar{s}b}, \
s\in{S},\ a,b\in{X^*}\}$ is a $k$-basis for $Lie(X|S)$.
\end{enumerate}
\end{theorem}
{\it Proof.} $(i)\Longrightarrow(ii)$. \ By noting that
$Id_{Lie}(S)\subseteq Id_{ass}(S)$, where $Id_{ass}(S)$ is the ideal
of $k\langle X\rangle$ generated by $S$, and by using Theorem 5.6
and Lemma 4.4, the result follows.

$(ii)\Longrightarrow(iii)$. \ Suppose that $\sum\limits_{[u_i]\in
Red(S)}\alpha_i[u_i]=0$ in $Lie(X|S)$ with $u_1>u_2>\cdots$, that
is, $\sum\limits_{[u_i]\in Red(S)}\alpha_i[u_i]\in{Id_{Lie}(S)}$.
Then each $\alpha_i$ must be 0. Otherwise, say $\alpha_1\neq0$.
Then,  by (ii), we know that
$\overline{\sum\limits_{i}\alpha_i[u_i]}=u_1$ which implies that
$[u_1]\not\in{Red(S)}$, a contradiction.

On the other hand, for any $f\in{Lie(X)}$, by Lemma 5.7, we have
$$
f+Id_{Lie}(S)=\sum\limits_i \alpha_i([u_i]+Id_{Lie}(S)).
$$

$(iii)\Longrightarrow(i)$. \ For any composition $\langle
f,g\rangle_w$ with $f,g\in{S}$, we have $\langle
f,g\rangle_w\in{Id_{Lie}(S)}$. Then, by (iii) and by Lemma 5.7,
$$
 \langle
f,g\rangle_w=\sum\beta_j[a_js_jb_j]_{\bar{s_j}},
$$
where each $\beta_j\in k, \ [a_js_jb_j]_{\bar{s_j}}$ is normal
$S$-word and $a_j\overline{s_j}b_j<w$. This proves that $S$ is a
Gr\"{o}bner-Shirshov basis in $Lie(X)$.

$(ii)\Longleftrightarrow(ii')$. \ This part is clear. \ \ \ \
$\square$


\begin{thebibliography}{99}
\bibitem{b}
G. M. Bergman: The diamond lemma for ring theory, {\it  Adv. in
Math.}{\bf\ 29} 178-218 (1978).

\bibitem{b72}
L. A. Bokut: Unsolvability of the word problem, and subalgebras of
finitely presented Lie algebras, {\it Izv. Akad. Nauk. SSSR Ser.
Mat.}{\bf\ 36} 1173-1219 (1972).

\bibitem{b76}
L. A. Bokut: Imbeddings into simple associative algebras, {\it
Algebra i Logika}{\bf\ 15} 117-142 (1976).

\bibitem{b99}L. A. Bokut, S.-J. Kang, K.-H. Lee, P. Malcolmson: Gr\"obner-Shirshov
bases for Lie super-algebras and their universal enveloping
algebras,  {\it J. Algebra}, {\bf 217 } 461-495 (1999).

\bibitem{b99-}L. A. Bokut, P. Malcolmson: Gr\"obner-Shirshov bases for
relations of a Lie algebra and its enveloping algebra, {\it Algebra
and Combinatorics } (Hong Kong), Springer, Singapore, 47-54 (1999).

\bibitem{bs}
 L. A. Bokut, K. P. Shum: Gr\"{o}bner and Gr\"{o}bner-Shirshov
bases in algebra: an elmentary approach,  {\it Southeast Asian
Bulletin of Mathematics}{\bf\ 29} 227-252 (2005).

\bibitem{bu65}
B. Buchberger: {\it An Algorithm for Finding a Basis for the Residue
Class Ring of a Zero-dimensional Polynomial Ideal}, Ph.D. thesis,
University of Innsbruck, Austria, (1965). (in German)

\bibitem{bu70}B. Buchberger: An algorithmical criteria for the
solvability of algebraic systems of equations, {\it Aequationes
Math.}, {\bf 4}  374-383 (1970). (in German)

\bibitem{C} K. T. Chen, R. H. Fox, and R. C. Lyndon: Free differential calculus, IV:
the quotient groups of the lower central series. {\it Annals of
Mathematics}, {\bf 68 } 81-95 (1958).


\bibitem{Cohn} P. M. Cohn:  review Zbl 0080.25503  Shirshov, A. I. \"Uber freie
Liesche Ringe.  {\it Mat. Sb.}, N. Ser.  {\bf 45 } (87) 113-122
(1958). (in Russian)

\bibitem{L} M. Lothaire:  ``Combinatorics on Words", Encyclopedia of
 Mathematics,  Vol. 17,
Addison-Wesley, 1983; reprented  M. Lothaire,  "Combinatorics on
Words", Cambridge Mathematical Library,  Cambridge University Press,
1997.


\bibitem{M} W. Magnus: \"Uber Beziehungen zwischern h\"oren
Kommutatoren, {\it J. Reine  Angew. Math,} {\bf 177} 105-115 (1937).

\bibitem{R} C. Reutenauer: {\it Free Lie Algebras},
  Oxford Science Publications, 1993.

\bibitem{S53}A. I. Shirshov: Subalgebras of free Lie algebras,
 {\it Mat. Sb.}, {\bf 33 } 441-452 (1953). (in Russian)


\bibitem{S58}
 A. I. Shirshov: On free Lie rings, {\it Sibirsk. Mat.
Sb.}{\bf\ 45}  (87) 113-122 (1958). (in Russian)


\bibitem{S62-}A. I. Shirshov: Bases for free Lie algebras, {\it Algebra
Logic,} {\bf 1} 14-19 (1962).


\bibitem{Sh}
A. I. Shirshov:  Some algorithmic problem for Lie algebras,
  {\it Sibirsk. Mat. Z.}{\bf\ 3}, 292-296 (1962) (in Russian);
   English translation in SIGSAM Bull.{\bf\ 33} (2) 3-6 (1999).

 \bibitem{W}E. Witt: Treue Darstellungen Lieschen Ringe,
 {\it J. Reine  Angew. Math.,} {\bf 177}  152-160 (1937).
\end{thebibliography}
\end{document}